 \documentstyle[graphicx,amscd,amssymb,verbatim,12pt,righttag,color]{amsart}
 \newtheorem{theorem}{Theorem}[section]
 
 \newtheorem{Prop}[theorem]{Proposition}
 \newtheorem{Lem}[theorem]{Lemma}
 \newtheorem{Cor}[theorem]{Corollary}
 
 \newtheorem{Example}[theorem]{Example}

 \numberwithin{equation}{section}
 \renewcommand{\rm}{\normalshape}
 
\date{}

\begin{document}
\title{Gabor orthonormal bases generated by the unit cubes}

\author{Jean-Pierre Gabardo}
\email{gabardo@@mcmaster.ca}

\address{Department of Mathematics and Statistics, McMaster University,
Hamilton, Ontario, L8S 4K1, Canada}

\author{Chun-Kit Lai}
 \email{cklai@@sfsu.edu}

\address{Department of Mathematics, San Francisco State University, 1600 Holloway Ave., San Francisco, CA 94132.
}

\author{Yang Wang}
 \email{yangwang@@ust.hk}

\address{Department of Mathematics, Hong Kong University of Science and Technology, Hong Kong}


\maketitle

\begin{abstract}
We consider the problem in determining the countable sets $\Lambda$ in the time-frequency plane
such that the Gabor system generated by the time-frequency shifts of the window $\chi_{[0,1]^d}$
 associated with $\Lambda$
forms a Gabor orthonormal basis for $ L^2({\Bbb R}^d)$.
We show that, if this is the case, the translates by elements $\Lambda$ of the unit cube
in ${\Bbb R}^{2d}$ must  tile  the time-frequency space ${\Bbb R}^{2d}$.
By studying the possible structure of such tiling sets,
we completely classify all such admissible sets $\Lambda$ of time-frequency shifts when $d=1,2$.
Moreover, an inductive procedure for constructing such sets $\Lambda$
in dimension $d\ge 3$ is also given. An interesting and surprising consequence of our
results is the existence, for $d\geq 2$,
of discrete sets $\Lambda$ with ${\mathcal G}(\chi_{[0,1]^d},\Lambda)$ forming a
Gabor orthonormal basis but with the associated ``time''-translates of the window
$\chi_{[0,1]^d}$ having significant overlaps.
\end{abstract}

\section{Introduction}

Let $g$ be a non-zero function in $L^2({\Bbb R}^d)$ and let $\Lambda$ be a discrete countable set
on ${\Bbb R}^{2d}$, where we identify ${\Bbb R}^{2d}$
to  the time-frequency plane by writing $(t,\lambda)\in\Lambda$ with $t,\lambda\in{\Bbb R}^d$.
 The  Gabor system associated with the window  $g$
consists of the set of translates and modulates of  $g$:
\begin{equation}\label{Gabor}
{\mathcal G}(g,\Lambda) = \{e^{2\pi i \langle\lambda,x\rangle}g(x-t): (t,\lambda)\in\Lambda\}.
\end{equation}
Such systems were first introduced by  Gabor \cite{[Gab]} who used them for applications in the theory of telecommunication,
but there has been a more recent interest in using Gabor system to  expand
functions both from a theoretical and applied perspective.
The branch of Fourier analysis dealing with Gabor systems is usually
referred  to  as Gabor, or time-frequency, analysis.
Gr\"ochenig's monograph \cite{[G]}
provide an excellent and detailed exposition on this subject.

\medskip

Recall that the Gabor system is a {\it frame} for $L^2({\Bbb R}^d)$  if there exists constants $A,B>0$ such that
 \begin{equation}\label{eq0.0}
A\|f\|^2\leq \sum_{(t,\lambda)\in\Lambda}|\langle f,
e^{2\pi i \langle\lambda,\cdot\rangle}g(\cdot-t)\rangle|^2\leq B\|f\|^2, \quad f\in L^2({\Bbb R}^d).
\end{equation}
It is called an orthonormal basis  for $L^2({\Bbb R}^d)$ if it is complete
and the elements of the system (\ref{Gabor}) are mutually orthogonal
in $L^2({\Bbb R}^d)$ and have norm 1, or,
equivalently, $\|g\|=1$ and  $A=B=1$ in (\ref{eq0.0}). One of the fundamental
problems in Gabor analysis is to classify the windows $g$ and
time-frequency sets $\Lambda$ with the property that the associated  Gabor
system ${\mathcal G}(g,\Lambda)$ forms a
(Gabor) frame or an orthonormal basis for $L^2({\Bbb R}^d)$.  This is of course
 a very difficult problem
and only partial results are known.
For example, to the best of our knowledge, the complete characterization of
time-frequency sets $\Lambda$ for which (\ref{Gabor}) is a frame
for $L^2({\Bbb R}^d)$  was only done
when $g = e^{-\pi x^2}$, the  Gaussian window. Lyubarskii, and Seip and
Wallsten \cite{[L],[SW]} showed
that ${\mathcal G}( e^{-\pi x^2},\Lambda)$ is a Gabor frame if and only
if the lower Beurling density of $\Lambda$ is strictly greater than  $1$.
If we assume that $\Lambda$ is a lattice of the form $a{\Bbb Z}\times b{\Bbb Z}$,
then it is well known that $a b\leq 1$ is a necessary condition for
(\ref{Gabor}) to form a frame for $L^2({\Bbb R}^d)$. Gr\"{o}chenig  and
St\"{o}ckler \cite{[GS]} showed that for totally positive
functions, (\ref{Gabor}) is a frame if and only if $ab<1$. If we consider
$g = \chi_{[0,c)}$, the characteristic function of an interval, the associated characterization problem
is known as the  {\it abc-problem} in Gabor analysis. By rescaling, one may assume that $c=1$.
In that case, the famous Janssen tie showed that the structure of the set of couples $(a,b)$ yielding a frame
is very complicated \cite{[J1],[GH]}. A complete solution of the abc-problem was recently obtained
by Dai and Sun \cite{[DS]}.

\medskip

In this paper, we focus our attention on
 Gabor system of the form (\ref{Gabor}) which yield orthonormal bases for $L^2({\Bbb R}^d)$.
 Perhaps the  most natural and simplest
example of Gabor orthonormal basis
is the system ${\mathcal G}(\chi_{[0,1]^d}, {\Bbb Z}^d\times{\Bbb Z}^d)$.
The orthonormality property for this system easily follows
 from that facts that the Euclidean space ${\Bbb R}^d$ can be partitioned by
the ${\Bbb Z}^d$-translates of the hypercube $[0,1]^d$ and that
the  exponentials $e^{2\pi i \langle n,x\rangle}$ form an orthonormal basis
for the space of square-integrable functions supported on
 any of these translated hypercubes.
A direct generalization of this observation is the following:

\begin{Prop}\label{prop0.1}
Let $|g| = |K|^{-1/2}\chi_{K}$, where $|\cdot|$ denotes the Lebesgue measure, and
$K\subset {\Bbb R}^d $ is measurable with finite Lebesgue measure.
Suppose that
\begin{itemize}
\item The translates of $K$ by the discrete set ${\mathcal J}$ are pairwise a.e.~disjoint
and cover ${\Bbb R}^d$ up to a set of zero measure.
\item For each $t\in{\mathcal J}$, the set of exponentials
$\{e^{2\pi i \langle\lambda,x\rangle}: \lambda\in\Lambda_t\}$
is an orthonormal basis for $L^2(K)$.
\end{itemize}
Let
\begin{equation}\label{Standard}
\Lambda = \bigcup_{t\in{\mathcal J}}\{t\}\times\Lambda_t.
\end{equation}
 Then ${\mathcal G}(g,\Lambda)$ is a Gabor orthonormal basis for $L^2({\Bbb R}^d)$.
\end{Prop}
Although its proof is straighforward and will be omitted (see also \cite{[LiW]}),
this proposition gives us a flexible way of constructing large families of
Gabor orthonormal basis.
The first condition above means that
$K$ is a {\it  translational tile} (with ${\mathcal J}$ called an associated {\it tiling set})
and the second one that $L^2(K)$ admits an  orthonormal basis of exponentials.
If this last condition holds, $K$ is called a {\it  spectral set}
(and each ${\Lambda}_t$ is an associated {\it spectrum}).
The connection between
translational tiles and spectral sets is quite mysterious.
They were in fact conjectured to be the same class of sets by Fuglede \cite{[Fu]},
but that statement was later disproved by Tao \cite{[T]}
and the exact relationship between the two classes remains unclear.

\medskip

For the fixed window $g_d = \chi_{[0,1]^d}$, we call a countable set
$\Lambda\subset{\Bbb R}^{2d}$ {\it standard} if it is  of the form
(\ref{Standard}). Motivated by the complete solution to the $abc$-problem,
our main objective in this paper is  to characterize the discrete sets $\Lambda$
(not necessarily lattices) with the property that the Gabor system
${\mathcal G}(g_d,\Lambda)$ is a Gabor orthonormal basis. First,
by generalizing the notion of {\it orthogonal packing region} (see Section 2) in the work of
Lagarias, Reeds and Wang \cite{[LRW]} to the setting of Gabor systems,
we deduce a general criterion for
${\mathcal G}(g_d,\Lambda)$ to be a Gabor orthonormal basis.

\begin{theorem}\label{th0.1}
${\mathcal G}(g_d,\Lambda)$ is a Gabor orthonormal basis if and only
if ${\mathcal G}(g_d,\Lambda)$ is an orthogonal set and
the translates of $[0,1]^d$ by the elements of $\Lambda$  tile  ${\Bbb R}^{2d}$.
\end{theorem}

\medskip

This criterion offers a very simple solution to our problem in the one-dimensional case.

\begin{theorem}\label{th0.2}
In dimension $d=1$,  the system ${\mathcal G}(g_1,\Lambda)$ is a Gabor orthonormal basis
if and only if $\Lambda$ is standard.
\end{theorem}

However, such a simple characterization ceases to exist in higher dimensions.
We will introduce an inductive procedure which allows us to construct a Gabor
orthonormal basis with window $g_d$ from  a Gabor orthonormal basis with window $g_n$, $n<d$.
This procedure can be used to produce many non-standard Gabor orthonormal basis and we call
a set $\Lambda$ obtained through this procedure {\it pseudo-standard}.
Assuming a mild condition on a low-dimensional time-frequency space,
we show that ${\mathcal G}(g_d,\Lambda)$ are essentially pseudo-standard
(See Theorem \ref{pseudo th}).

Although we do not have a complete description of the sets $\Lambda$ yielding
Gabor orthonormal bases with window $g_d$ in dimension $d\ge 3$,
 we managed to obtain a complete characterization of those discrete sets $\Lambda\subset {\Bbb R}^4$
such that ${\mathcal G}(g_2,\Lambda)$ form an orthonormal basis for $L^2({\Bbb R}^2)$.

\begin{theorem}\label{th0.3}
${\mathcal G}(\chi_{[0,1]^2},\Lambda)$ is a  Gabor orthonormal basis for
 $L^2({\Bbb R}^2)$ if and only if we can partition ${\Bbb Z}$ into ${\mathcal J}$ and ${\mathcal J}'$ such that either
$$
\Lambda=\bigcup_{n\in{\mathcal J}} \{(m+t_{n,k},n, j+\mu_{k,m,n}, k+\nu_n): m,j,k\in{\Bbb Z}\}\cup \bigcup_{m\in{\Bbb Z}}\bigcup_{n\in{\mathcal J}'}\{(m+t_n,n)\}\times\Lambda_{m,n}
$$
or
$$
\Lambda= \bigcup_{m\in{\mathcal J}}\{(m,n+t_{m,j}, j+\nu_m, k+\mu_{j,m,n}): n,j,k\in{\Bbb Z}\}\cup \bigcup_{n\in{\Bbb Z}}\bigcup_{m\in{\mathcal J}'}\{(m,n+t_m)\}\times\Lambda_{m,n}.
$$
where $\Lambda_{m,n}+[0,1]^2$ tile ${\Bbb R}^2$ and $t_{n,k}$, $\mu_{k,m,n}$ and $\nu_n$ are real numbers in $[0,1)$ as a function of $m,n$ or $k$.
\end{theorem}

\medskip

\begin{figure}[h]
  \includegraphics[width=4in]{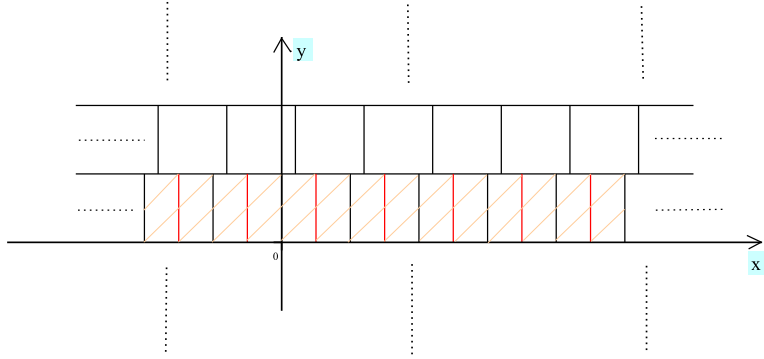}\\
  \caption{This figure illustrates the time-domain of $\Lambda$ in the first situation of Theorem \ref{th0.3}. We basically partition ${\Bbb R}^2$ by horizontal strips. Some strips, like ${\Bbb R}\times[0,1]$ with $n=0$, have overlapping structure. This corresponds to the first union of $\Lambda$. Some strips, like ${\Bbb R}\times[1,2]$ with $n=1$, have tiling structures. This corresponds to the second union of $\Lambda$.  }
\end{figure}

\medskip

We organize the paper as follows. In Section 2, we provide some preliminaries notations and
prove Theorem \ref{th0.1}. In Section 3, we prove Theorem \ref{th0.2} and introduce the
pseudo-standard time-frequency set. In the last section, we focus on dimension 2 and prove Theorem \ref{th0.3}.

%
%
%
%

\bigskip

\section{Preliminaries}
In this section, we explore the relationship between Gabor orthonormal bases and
tilings in the time-frequency space. This theory will be an extension of spectral-tile duality in \cite{[LRW]} to the setting of Gabor analysis.
Denote by $|K|$ the Lebesgue measure of a set $K$. We say that a closed set $T$ is a {\it region}
 if $|\partial T|=0$ and $\overline{T^{o}} = T$. A bounded region $T$ is called a
 {\it translational tile} if we can find a countable set ${\mathcal J}$ such that
\begin{enumerate}
\item $|(T+t)\cap(T+t')|=0,\quad t,t'\in \mathcal{J}$, $t\ne t'$, and
\item $\bigcup_{t\in{\mathcal J}}(T+t) = {\Bbb R}^d$.
\end{enumerate}
In that case,  ${\mathcal J}$ is called a {\it tiling set} for $T$ and
 $T+{\mathcal J}$ a tiling of ${\Bbb R}^d$. We will say that $T+{\mathcal J}$ is a packing of
${\Bbb R}^n$ if (1) above is satisfied.
  We can generalize the notion of tiling and packing to measures and functions.
 Given a positive Borel measure $\mu$ and $f\in L^1({\Bbb R}^n)$ with $f\ge 0$,
 the convolution of $f$ and $\mu$ is defined to be
$$
f\ast\mu(x) = \int\, f(x-y)\,d\mu(y),\quad x\in {\Bbb R}^n,
$$
(where a Borel measurable function is chosen in the equivalence class of $f$ to define the integral
above).
We say that {\it $f+\mu$ is a tiling (resp. packing) of
${\Bbb R}^d$} if $f\ast \mu =1$ (resp. $f\ast \mu\leq 1$) almost everywhere with respect to the
Lebesgue measure. It is clear that if $f = \chi_T$ and
 $\mu = \delta_{{\mathcal J}}$ where $\delta_{{\mathcal J}}=\sum_{t\in\mathcal J}\,\delta_t$,
 then $f\ast\mu = 1$ is equivalent to $T+{\mathcal J}$ being a tiling.

\medskip

First, we start with the following theorem which gives us a very useful
 criterion to decide if a packing is actually a tiling. In fact,
special cases of this theorem were proved  by many different authors in different settings
(see e.g. \cite[Theorem 3.1]{[LRW]}, \cite[Lemma 3.1]{[K]} and  \cite{[Li]}), but
the following version is the most general one as far as we know.
\begin{theorem}\label{tiling Prop}
Suppose that $F,G \in L^1({\mathbb{R}}^n)$ are two functions with $F,G\ge 0$
and $\int_{{\mathbb{R}}^n}\,F(x)\,dx=\int_{{\mathbb{R}}^n}\,G(x)\,dx=1$.
Suppose that $\mu$ is a positive Borel measure
on ${\mathbb{R}}^n$ such that
$$
F*\mu\le 1 \quad \text{and}\quad G*\mu\le 1.
$$
Then,  $F*\mu=1$ if and only if $G*\mu=1$.
\end{theorem}

\begin{pf}
By symmetry, it suffices to prove one side of the equivalence.
Assuming that $F*\mu=1$, we have
$$
1=F*\mu \ \Rightarrow  \ 1 = 1*G = G*F*\mu=F*G*\mu.
$$
Letting $H=G*\mu$ we have $0\le H\le 1$ and $H*F=1$. We now show that $H=1$. Indeed letting $A$
be the set $\{x\in{\mathbb{R}}^n, H(x)<1\}$ and $B={\mathbb{R}}^n\setminus A$,
we have
$$
(H*F)(x)=\int_{{\mathbb{R}}^n}\,H(y)\,F(x-y)\,dy
=\int_{A}\,H(y)\,F(x-y)\,dy+\int_{B}\,H(y)\,F(x-y)\,dy
$$
Now, if $|A|>0$, we have
$$
\int_{{\mathbb{R}}^n}\int_{A}\,F(x-y)\,dy\,dx= |A|>0
$$
and there exists thus a set $E$ with positive measure such that
$$
\int_{A}\,F(x-y)\,dy>0,\quad x\in E.
$$
If $x\in E$, we have
$$\begin{aligned}
 \int_{A}\,H(y)\,F(x-y)\,dy+\int_{B}\,H(y)\,F(x-y)\,dy
<&\int_{A}\,F(x-y)\,dy+\int_{B}\,F(x-y)\,dy\\
=&(1*F)(x)=1.\\
\end{aligned}$$
This contradicts to the fact that  $H*F=1$ almost everywhere. Hence, $|A|=0$ and $H=1$ follows.
\end{pf}

\medskip

Let $f,g\in L^2({\Bbb R}^d)$. We define the {\it short time Fourier transform} of $f$ with respect
to the window $g$ be
$$
V_gf(t,\nu) = \int_{{\Bbb R}^{2d}}\, f(x)\,\overline{g(x-t)}\,e^{-2\pi i\langle\nu,x\rangle}\,dx.
$$
Let ${\mathcal G}(g,\Lambda)$ be a Gabor orthonormal basis.
Since translating $\Lambda$ be an element of  ${\Bbb R}^{2d}$ does not affect the
orthonormality nor the completeness of the given system,
there is no loss of generality in assuming that $(0,0)\in\Lambda$.
We say that a region $D$ ($\subset{\Bbb R}^{2d}$) is an {\it orthogonal packing region} for $g$ if
$$
(D^{\circ}-D^{\circ}) \cap {\mathcal Z}(V_gg) = \emptyset.
$$
Here ${\mathcal Z}(V_gg) = \{(t,\nu): V_gg(t,\nu)=0\}$.
\medskip

\medskip

\begin{Lem}\label{lem2.1.1}
Suppose that ${\mathcal G}(g,\Lambda)$ is a mutually orthogonal set of $L^2({\Bbb R}^d)$. Let $D$ be any
orthogonal packing region for $g$. Then $\Lambda-\Lambda\subset {\mathcal Z}(V_gg)\cup\{0\}$ and
$\Lambda+D$ is a packing of ${\Bbb R}^{2d}$. Suppose furthermore that ${\mathcal G}(g,\Lambda)$
is a Gabor orthonormal basis. Then $|D|\leq 1$.
\end{Lem}

\medskip

\begin{pf}
Let $(t,\lambda),(t',\lambda')\in\Lambda$ be two distinct points in $\Lambda$. Then
 $$
 \int\, g(x-t')\,\overline{g(x-t)}\,e^{-2\pi i (\lambda-\lambda')x}\,dx = 0,
 $$
 or equivalently,  after the change of variable $y=x-t'$,
 $$
\int\, g(x)\,\overline{g(x-(t-t'))}\,e^{-2\pi i (\lambda-\lambda')x}\,dx = 0.
$$
Hence, $V_gg(t-t',\lambda-\lambda')=0$ and $(t,\lambda)-(t',\lambda')\in{\mathcal Z}(V_gg)$.
This means that  $(t,\lambda)-(t',\lambda')\not\in D^\circ-D^\circ$. Therefore, the intersection
of the sets
$(t,\lambda)+D$  and $(t',\lambda')+D$ has zero Lebesgue measure.

\medskip

Suppose now that  ${\mathcal G}(g,\Lambda)$ is a Gabor orthonormal basis.
Denote by $R$ the diameter of $D$. By the packing property of $\Lambda+D$,
$$
\begin{aligned}
|D|\cdot\frac{\#(\Lambda\cap [-T,T]^{2d})}{(2T)^{2d}}
 =& \frac{1}{(2T)^{2d}}\left|\bigcup_{\lambda\in\Lambda\cap[-T,T]^{2d}}(D+\lambda)\right|\\
\leq&\frac{1}{(2T)^{2d}} \left|[-T-R,T+R]^{2d}\right| = (1+\frac{R}{T})^{2d}.
\end{aligned}
$$
Taking limit $T\rightarrow\infty$ and using the fact that
Beurling density of $\Lambda$ is 1 (\cite{[RS]}), we have $|D|\leq 1$.
\end{pf}

\medskip

We say that an orthogonal packing region $D$ for $g$
is {\it tight} if  we have furthermore $|D|=1$.
We now apply Theorem \ref{tiling Prop} to the Gabor
orthonormal basis problem.

\begin{theorem}\label{th2.1}
Suppose that ${\mathcal G}(g,\Lambda)$ is an orthonormal
set in $L^2({\Bbb R}^d)$ and that
$D$ is a tight orthogonal packing region for $g$.
Then ${\mathcal G}(g,\Lambda)$ is a Gabor orthonormal basis for $L^2({\Bbb R}^d)$
if and only if $\Lambda+D$ is a tiling of ${\Bbb R}^{2d}$.
\end{theorem}

\begin{pf}
Let $F = \chi_{D}$ and $G = |V_gf|^2/\|f\|_2^2$.
Then $\int_{{\Bbb R}^{2d}}F = 1$ and $\int_{{\Bbb R}^{2d}}G = \|g\|_2^2=1$.
 Now, as $D$ is an
orthogonal packing region for $g$,  we have in particular
$$
\sum_{\lambda\in{\Lambda}}\chi_{D}(x-\lambda)\leq 1.
$$
This shows that
$$
\delta_{\Lambda}\ast F = \delta_{\Lambda}\ast \chi_D \leq 1.
$$
Moreover, $\Lambda+D$ is a tiling of  ${\Bbb R}^{2d}$ if and only
if $\delta_{\Lambda}\ast \chi_D = 1$.
 On the other hand, $(g,\Lambda)$ being  a mutually orthogonal set,
Bessel's inequality  yields
$$
\sum_{(t,\lambda)\in\Lambda}\,
\left|\int_{{\Bbb R}^d}\, f(x)\,\overline{g(x-t)}\,
e^{-2\pi i \langle\lambda,x\rangle}\,dx\right|^2
\leq \|f\|^2,\quad f\in L^2({\Bbb R}^d),
$$
or, replacing $f$ by $f(x-\tau)e^{2\pi i \nu x}$ with $(\tau,\nu)\in {\Bbb R}^{2d}$,
$$
\sum_{(t,\lambda)\in\Lambda}\left|V_gf(\tau-t,\nu-\lambda)\right|^2
\leq \|f\|^2,\quad f\in L^2({\Bbb R}^d).
$$
Hence,
$$
\delta_{\Lambda}\ast G = \delta_{\Lambda}\ast \frac{|V_gf|^2}{\|f\|^2}\leq 1
$$
with equality  if and only if the Gabor orthonormal system is in fact a basis.
The conclusion follows then from Theorem \ref{tiling Prop}.
\end{pf}

\medskip

\noindent{\bf Proof of Theorem \ref{th0.1}.}  Let $g_d = \chi_{[0,1]^d}$.
Using Theorem \ref{th2.1},
we just need to show that
$[0,1]^{2d}$ is a
tight orthogonal packing region for $g_d$.

\medskip

We first consider the case $d=1$. For $g_1 = \chi_{[0,1]}$,  a
direct computation shows that
\begin{equation}\label{Vgg}
V_{g_1}{g_1}(t,\nu) = \left\{
              \begin{array}{ll}
                0, & \hbox{$|t|\geq1$;} \\
                \frac{1}{2\pi i \nu}\left(e^{2\pi i \nu t}-e^{2\pi i \nu}\right),
 & \hbox{$0\leq t\leq 1$;} \\
                \frac{1}{2\pi i \nu}\left(1-e^{2\pi i \nu(t+1)}\right),
& \hbox{$-1\leq t \leq 0$.}
              \end{array}
            \right.
\end{equation}
The zero set of $V_{g_1}g_1)$ is therefore given by
\begin{equation}\label{ZVgg}
{\mathcal Z}(V_{g_1}g_1) = \{(t,\nu): |t|\geq 1\}\cup \{(t,\nu):
 \nu(1-|t|)\in {\Bbb Z}\setminus \{0\}\}.
\end{equation}
Hence, $(0,1)^2-(0,1)^2 = (-1,1)^{2}$ does not intersect the zero set and therefore
$[0,1]^2$ is a tight orthogonal packing region for $g_1$.

\medskip

We now consider the case $d\ge 2$. As we can decompose $g_d$ as
$\chi_{[0,1]}(x_1)...\chi_{[0,1]}(x_d)$, we have
$$
V_{g_d}g_d(t,\nu) = V_{g_1}g_1(t_1,\nu_1)\dots
V_{g_1}g_1(t_d,\nu_d)\,\,\text{where}\,\,\,t = (t_1,\dots,t_d)
\,\,\text{and}\,\,\nu = (\nu_1,\dots,\nu_d).
$$
The  zero set of $V_{g_d}g_d$ is therefore given by
\begin{equation}\label{ZVggd}
{\mathcal Z}(V_{g_d}g_d) = \{(t,v): |t|_{\max}\geq 1\}\cup
\left(\bigcup_{i=1}^{d}\{(t,\nu): \nu_i(1-|t_i|)\in{\Bbb Z}\setminus\{0\})\}\right)
\end{equation}
where $|t|_{\max} = \max\{t_1,...,t_d\}$. It follows that $[0,1]^{2d}$ is a tight
orthogonal packing region for $g_d$.
\qquad$\Box$

\medskip

The following example will not be used in later discussion, but it
 demonstrates the usefulness of the theory for windows other than the unit cube.
\begin{Example}
{\rm Let $g(x) = \frac{2}{e^{2x}+e^{-2x}}$ be the hyperbolic secant function.
It can be shown (\cite{[J2]}; see also \cite{[Ga]}) that}
$$
V_gg(t,\nu) = \frac{\pi\sin(\pi \nu t)e^{-\pi i \nu t}}{\sinh(2t)\sinh(\pi^2\nu/2)}
$$
{\rm and the zero set is given by}
$$
{\mathcal Z}(V_gg) = \{(t,\nu): t\nu\in{\Bbb Z}\setminus\{0\}\}.
$$
{\rm Hence, $[0,1]^2$ is a tight orthogonal packing region for $g$.
Note that the zero set does not contain any points on the $x$- axis and $y$-axis.
There is no tiling set $\Lambda$ for  $[0,1]^2$
such that $\Lambda-\Lambda\subset {\mathcal Z}(V_gg) \cup\{0\}$ (see also Proposition \ref{prop3.2}
in the next section) and thus there is no Gabor
orthonormal basis using the hyperbolic secant as a window. This can be viewed as a particular
case of a version of the Balian-Low theorem valid for irregular Gabor frames
which was recently obtained in \cite{[AFK]} and which state that Gabor orthonormal bases cannot exist
if the window function is in the modulation space $M^1({\Bbb R}^d)$.}
\end{Example}

\medskip

\section{Gabor orthonormal bases}

Using Lemma \ref{lem2.1.1}, Theorem \ref{th0.1} may be restated in the following way:

\begin{theorem}\label{prop3.1}
 ${\mathcal G}(\chi_{[0,1]^d},\Lambda)$ is a Gabor orthonormal basis if and only if
the inclusion $\Lambda-\Lambda\subset{\mathcal Z}(V_gg)\cup\{0\}$ holds and $\Lambda+[0,1]^{2d}$ is a tiling.
\end{theorem}

\medskip

In view of the previous result, the possible translational tilings of the unit cube on ${\Bbb R}^{2d}$
play a fundamental role in the solution of our problem.
A characterization for these is not available in arbitrary $2d$ dimension
but it is easily obtained when $d=1$. We prove this result here for completeness
 but it should be well known.

\begin{Prop}\label{prop3.2}
Suppose that $\chi_{[0,1]^2}+{\mathcal J}$ is a  tiling  of ${\Bbb R}^2$ with $(0,0)\in{\mathcal J}$.
Then ${\mathcal J}$ is of either of the following two form:
\begin{equation}\label{2Dtilingset}
{\mathcal J} = \bigcup_{k\in{\Bbb Z}} ({\Bbb Z}+a_k)\times\{k\} \ \mbox{or}  \ {\mathcal J}
= \bigcup_{k\in{\Bbb Z}}\{k\}\times({\Bbb Z}+a_k)
\end{equation}
where $a_k$ are any real numbers in $[0,1)$ for $k\ne 0$ and $a_0=0$.
\end{Prop}

\begin{pf}
By Keller's criterion  for square tilings (see e.~g. \cite[Proposition 4.1]{[LRW]}),
for any $(t_1,t_2)$ and $(t_1',t_2')$ in ${\mathcal J}$,
$t_i-t_i'\in{\Bbb Z}\setminus\{0\}$ for some $i=1,2$.
 Taking $(t_1',t_2')=(0,0)$, we obtain that, for any
$(t_1,t_2)\in {\mathcal J}\setminus \{(0,0)\}$,
one of $t_1$ or $t_2$
belongs to ${\Bbb Z}\setminus\{0\}.$ If ${\mathcal J} \subset {\Bbb Z}$, we must have
 ${\mathcal J}={\Bbb Z}$ for  $\chi_{[0,1]^2}+{\mathcal J}$ to be  tiling  of ${\Bbb R}^2$
and $ {\Bbb Z}$ can be written as
either of the sets in (\ref{2Dtilingset}) by taking $a_k=0$ for all $k$.
 Suppose that there exists $(s_1, s_2)\in{\mathcal J}$
such that $s_1$ is not an integer and $s_2\in \mathbb{Z}$.
If $(t_1,t_2)\in {\mathcal J}$
and $t_2 \notin {\Bbb Z}$, then both $t_1$ and $t_1-s_1$
must be integers which would imply
that $s_1$ is an integer, contrary to our assumption.
Hence, $(s_1, s_2)\in{\mathcal J}$ implies
$s_2\in {\Bbb Z}$ and  we can write
$$
{\mathcal J} = \bigcup_{k\in{\Bbb Z}} {\mathcal J}_k\times\{k\}.
$$
for some discrete set ${\mathcal J}_k\subset {\Bbb R}$. For
$\chi_{[0,1]^2}+{\mathcal J}$ to be a tiling  of ${\Bbb R}^2$,
the set ${\mathcal J}_k$
must be of the form ${\mathcal J}_k = {\Bbb Z}+a_k.$
In that case ${\mathcal J}$ can be expressed as one of the
sets in the first collection appearing
in (\ref{2Dtilingset}).

Similarly,  if there exists $(s_1, s_2)\in{\mathcal J}$
such that $s_2$ is not an integer and $s_1\in \mathbb{Z}$,
 ${\mathcal J}$ can be expressed as one of the sets in the
second collection appearing
in (\ref{2Dtilingset}).
This completes the proof.
\end{pf}

\medskip

We say that the Gabor orthonormal basis ${\mathcal G}(\chi_{[0,1]^d},\Lambda)$
is {\it standard} if
$$
\Lambda = \bigcup_{t\in{\mathcal J}}\{t\}\times\Lambda_t,
$$
where ${\mathcal J}+[0,1]^d$ tiles ${\Bbb R}^d$ and
$\Lambda_{t}$ is a spectrum for  $[0,1]^d$.
(Note that, by the result in \cite{[LRW]},
$\Lambda_t+[0,1]^d$ must then be a tiling of ${\Bbb R}^d$
for every $t\in{\mathcal J}$.)

The following result settles the one-dimensional case.

\begin{theorem}
 ${\mathcal G}(\chi_{[0,1]},\Lambda)$ is a Gabor orthonormal
 basis if and only if  $\Lambda$ is standard.
\end{theorem}

\begin{pf}
We just need to show that $\Lambda$ being standard is a necessary condition
for ${\mathcal G}(\chi_{[0,1]},\Lambda)$ to be a Gabor orthonormal basis.
We can also assume, for simplicity, that $(0,0)\in \Lambda$.
By Proposition \ref{prop3.1}, if ${\mathcal G}(\chi_{[0,1]},\Lambda)$
is a Gabor orthonormal basis,
then
$\Lambda-\Lambda\subset{\mathcal Z}(V_gg)\cup\{0\}$ and
$\Lambda+[0,1]^2$ must be  a tiling of ${\Bbb R}^2$. By Proposition \ref{prop3.2},
$\Lambda$ must be of either one of the forms  in (\ref{2Dtilingset}).
Note that  $\Lambda$ is standard in the second case.
In order to deal with the first case, suppose that
$$
\Lambda=\bigcup_{k\in{\Bbb Z}}\,({\Bbb Z}+a_k)\times\{k\},\,\,\text{with}\,\,
a_k\in[0,1),\,\,k\ne 0,\,\,a_0=0.
$$
We now show that this  is impossible unless $a_k = 0$ for all $k$
(which reduces to the case $\Lambda={\Bbb Z}^2$, which is standard).
We can assume, without loss of generality, that  $a_k\neq 0$ for some $k>0$
with $k$ being the smallest
such index.
If $a_k\ne 0$ for some $k$, then both $(a_k,k)$ and $(0,k-1)$ are in $\Lambda$.
The orthogonality of the Gabor system
then implies that $(a_k,1)\in{\mathcal Z}(V_gg)$. Using (\ref{ZVgg}),
we deduce that $1\cdot(1-|a_k|)\in{\Bbb Z}\setminus\{0\}.$
That means $a_k$ must be an integer, which is a contradiction.
 Hence, the first case is impossible unless $a_k=0$ for all $k$
 and the proof is completed.
\end{pf}

\medskip

A description of all time-frequency sets $\Lambda$ for
which ${\mathcal G}(\chi_{[0,1]^d},\Lambda)$ is a
Gabor orthonormal basis however
become vastly more complicated when $d\ge 2$.
In particular, as we will see, the standard structure
cannot cover all possible cases.
Consider integers $m,n >0$ such that $m+n=d$.
For convenience and to be consistent with our previous notation, we will write the cartesian product
of the two  time-frequency spaces ${\Bbb R}^{2m}$ and ${\Bbb R}^{2n}$ in the non-standard form
$$
{\Bbb R}^{2d}={\Bbb R}^{2m}\times{\Bbb R}^{2n}=\{(s,t,\lambda,\nu),\,\,(s,\lambda)\in {\Bbb R}^{2m},\,\,
(t,\nu)\in {\Bbb R}^{2n}\}.
$$
We will also denote by $\Pi_1$ the projection operator from ${\Bbb R}^{2d}$ to ${\Bbb R}^{2m}$
defined by
\begin{equation}\label{pidef}
\Pi_1\left((s,t,\lambda,\nu)\right)=(s,\lambda),\quad (s,t,\lambda,\nu)\in {\Bbb R}^{2d}
={\Bbb R}^{2m}\times{\Bbb R}^{2n}.
\end{equation}
To simplify the notation, we also define $g_k = \chi_{[0,1]^k}$
for  any $k\ge 1$.
We now build a new family of time-frequency sets on ${\Bbb R}^{2d}$ as follows.
Suppose that
${\mathcal G}(\chi_{[0,1]^{m}},\Lambda_1)$ is a Gabor  orthonormal basis
for $L^2({\Bbb R}^{m})$ and that we associate with each $(s,\lambda)\in \Lambda_1$,
a discrete set $\Lambda_{(s,\lambda)}$ in ${\Bbb R}^{2n}$ such that
 ${\mathcal G}(\chi_{[0,1]^{n}},\Lambda_{(s,\lambda)})$
is a Gabor orthonormal basis of $L^2({\Bbb R}^{n})$. We then define
\begin{equation}\label{eqLambda}
      \Lambda = \bigcup_{(s,\lambda)\in\Lambda_1}\,
\{(s,t,\lambda,\nu),\,\,(t,\nu)\in \Lambda_{(s,\lambda)}\}.
\end{equation}
 We say that a Gabor system  ${\mathcal G}(\chi_{[0,1]^{d}},\Lambda)$
with $\Lambda$ as in (\ref{eqLambda}) is  {\it pseudo-standard}.

\medskip

\begin{Prop}\label{prop3.4}
Every pseudo-standard Gabor system ${\mathcal G}(\chi_{[0,1]^d},\Lambda)$ is a
Gabor orthonormal basis of $L^2({\Bbb R}^d)$.
\end{Prop}

\medskip

\begin{pf}
If $ x\in{\Bbb R}^{m}$ and $y\in {\Bbb R}^{n}$, we have
$g_d (x,y) = g_m(x)g_n(y)$
(for $m+n=d$). This yields immediately that
\begin{equation}\label{tensor}
V_{g_d}g_d(s,t,\lambda,\nu)
= V_{g_{m}}g_{m}(s,
\lambda)\,V_{g_{n}}g_{n}(t,\nu),\quad
(s,\lambda)\in{\Bbb R}^{2m},
\,\,(t,\nu)\in{\Bbb R}^{2n}.
\end{equation}
Suppose that $\rho=(s,t, \lambda,\nu)$ and $\rho'=(s',t', \lambda',\nu')$
are distinct elements of $\Lambda$.
If  $(s, \lambda)=(s', \lambda')$,  then
 $(t, \nu)$ and $(t', \nu')$ are distinct elements of
$\Lambda_{(s,\lambda)}$ and we have thus
$$
(t'-t,\nu'- \nu)\in  {\mathcal Z}(V_{g_{n}}g_{n})
$$
which implies that ${\mathcal Z}(V_{g_{d}}g_{d}) (\rho'-\rho)=0$.
On the other hand, if  $(s, \lambda)\ne(s',\lambda'))$, we have then
$$
(s'-s,\lambda'- \lambda)\in  {\mathcal Z}(V_{g_{m}}g_{m})
$$
which implies again that ${\mathcal Z}(V_{g_{d}}g_{d}) (\rho'-\rho)=0$.
This proves the orthonormality of the system ${\mathcal G}(\chi_{[0,1]^d},\Lambda)$.
This proposition can now  be proved by invoking Theorem \ref{prop3.1} if we can show
that $\Lambda+[0,1]^{2d}$ is a tiling of ${\Bbb R}^{2d}$.
To prove this, we note that $\Lambda_1+[0,1]^{2m}$ is a
tiling of the subspace ${\Bbb R}^{2m}$
by Theorem \ref{prop3.1} and that, similarly, for each $(t,\lambda)\in $,
 $\Lambda_{(t,\lambda)}+[0,1]^{2n}$ is a tiling of ${\Bbb R}^{2n}$.
This easily implies the required tiling property and concludes the proof.
\end{pf}

\medskip

\begin{Example}
{\rm Consider the two-dimensional case $d=2$.
Let
$$
\Lambda_1 = \bigcup_{m\in{\Bbb Z}}\{m\}\times({\Bbb Z}+\mu_{m}),\quad \mu_m\in [0,1).
$$
Associate with each $\gamma = (m,j+\mu_m)\in\Lambda_1 $, the set
$$
\Lambda_{\gamma} = \bigcup_{n\in{\Bbb Z}}\{n+s_{m,j}\}\times({\Bbb Z}+\nu_{n,m,j}),
\quad s_{m,j}\in {\Bbb R} , \nu_{n,m,j}\in [0,1).
$$
Then,
$$
\Lambda:=  \{(m,n+s_{m,j},j+ \mu_m, k+\nu_{n,m,j}): m,n,j,k\in{\Bbb Z}\}
$$
(written in the form of $(t_1,t_2,\lambda_1,\lambda_2)$ where
$(t_1,t_2)$ are the translations and $(\lambda_1,\lambda_2)$ the frequencies) has the
pseudo-standard structure.
Note that the parameters  $s_{m,j}$ can be chosen so that the set $\Lambda$ is not standard
as the set
$$
\{(m,n+s_{m,j}),\,\,m,n,j\in{\Bbb Z}\}+[0,1]^2
$$
will not tile ${\Bbb R}^2$ in general. For example, for $m=n=0$, we could let
$s_{0,0}=0$ and
the numbers $s_{0,j}$
could be chosen as distinct numbers in the interval
$[0,1)$. The square $[0,1]^2$ would then overlap with
infinitely many of its translates appearing
as part of the Gabor system. }
\end{Example}

Using a similar procedure to higher dimension, we can produce many
non-standard Gabor orthonormal bases with window $\chi_{[0,1]^d}$.
However, the pseudo-standard structure still cannot cover all possible
cases of time-frequency sets. A time-frequency set could be a mixture of
pseudo-standard and standard structure. For example, consider the set
$$
\Lambda=\bigcup_{n\in{\Bbb Z}\setminus\{1\}}
\{(m+t_{n,k},n, j+\mu_{k,m,n}, k+\nu_n): j,k\in{\Bbb Z}\}\cup
\{(m,1)\}\times\Lambda_{m},
$$
where $\Lambda_m+[0,1]^2$ tiles ${\Bbb R}^2$. This set consists of two parts.
The first part is a subset of a set having the pseudo-standard
 structure while the second part
is a  subset of a set having the standard one.
Moreover, the translates of the unit square associated with the first part
are disjoint with those associated with the second part, showing that
${\mathcal G}(\chi_{[0,1]^2},\Lambda)$ is a mutually orthogonal set.
Since $\Lambda$  is clearly a tiling of ${\Bbb R}^4$, Theorem \ref{prop3.1}
shows that ${\mathcal G}(\chi_{[0,1]^2},\Lambda)$ is a Gabor orthonormal basis.

In the next section, we will classify all possible sets
$\Lambda\subset{\Bbb R}^4$
with the property that
${\mathcal G}(\chi_{[0,1]^2},\Lambda)$ is a Gabor orthonormal basis
for $L^2({\Bbb R}^2)$. However, we have

\begin{theorem}\label{pseudo th} Let $d=m+n$ and let $\Pi_1:{\Bbb R}^{2d}\to{\Bbb R}^{2m}$
be defined  by (\ref{pidef}).
Suppose that $(\chi_{[0,1]^d},\Lambda)$ is a Gabor orthonormal basis
and that $\Pi_1(\Lambda)+[0,1]^{2m}$ tiles ${\Bbb R}^{2m}$.
Then $\Lambda$ has the  pseudo-standard structure.
\end{theorem}

\begin{Prop}\label{prop3.3} Let $d=m+n$ and suppose that $(\chi_{[0,1]^d},\Lambda)$
is a Gabor orthonormal basis for $L^2({\Bbb R}^{d})$.
If $(s_0,\lambda_0)\in {\Bbb R}^{2m}$,
consider the translate of the  unit hypercube  in ${\Bbb R}^{2m}$,
 $C = (s_0,\lambda_0)+[0,1)^{2m}$, and
define
$$
\Lambda(C): = \{(t,\nu)\in{\Bbb R}^{2n}: (s,t,\lambda,\nu)
\in \Lambda \ \mbox{and} \ (s,\lambda)\in C\}.
$$
Then $(\chi_{[0,1]^n},\Lambda(C))$ is a Gabor orthonormal basis for $L^2({\Bbb R}^{2n})$.
\end{Prop}

\medskip

\begin{pf}
We first show that the system $(\chi_{[0,1]^n},\Lambda(C))$ is orthogonal.
Let $(t,\nu)$ and $(t',\nu')$ be distinct elements
of $\Lambda(C)$. There exist $(s,\lambda)$ and $(s',\lambda')$ in ${\Bbb R}^{2m}$
such that   $(s, t,\lambda,\nu)$ and  $(s', t',\lambda',\nu')$
both belong to $\Lambda$.
Using  the mutual
orthogonality of  the system $(\chi_{[0,1]^d},\Lambda)$ together with (\ref{tensor}), we have
$$
V_{g_m}g_m(s-s',\lambda-\lambda')=0 \ \mbox{or}
\ V_{g_n}g_n(t-t',\nu-\nu')=0.
$$
Note that, as both  $(s,\lambda)$ and  $(s',\lambda')$ belong
to $C$, we have
 $|s-s'|_{\max}<1$ and $|\lambda-\lambda'|_{\max}<1$.
In particular, $V_{g_m}g_m(s-s',\lambda-\lambda')\neq 0$
and the orthogonality of  the system $(\chi_{[0,1]^n},\Lambda(C))$
follows.

\medskip

If $(s,\lambda)\in \Pi_1(\Lambda)$ (as defined in (\ref{pidef})), let
$$
\Lambda_{(s,\lambda)} = \{(t,\nu):(s,t,\lambda,\nu)\in \Lambda \}.
$$
Let $f_1\in L^2({\Bbb R}^m)$, $f_2\in L^2({\Bbb R}^n)$ and
$(s_0,\lambda_0)\in {\Bbb R}^{2m}$.
Applying Parseval's identity to the function
$$
f(x,y)= e^{2\pi i \lambda_0\cdot x }\,f_1(x-s_0)\,f_2(y),\quad x\in {\Bbb R}^m,\,\,y\in{\Bbb R}^n,
$$
we obtain that
\begin{align*}
&\int_{{\Bbb R}^m}\,|f_1(x)|^2\,dx\,\int_{{\Bbb R}^n}\,|f_2(y)|^2\,dy\\
&= \sum_{(s,\lambda)\in\Pi_1(\Lambda)}
\sum_{(t,\nu)\in \Lambda_{(s,\lambda)}}
|V_{g_m}f_1(s-s_0,\lambda-\lambda_0)|^2\,|V_{g_n}f_2(t,\nu)|^2\\
&= \sum_{(s,\lambda)\in\Pi_1(\Lambda)}
\sum_{(t,\nu)\in \Lambda_{(s,\lambda)}}
|V_{f_1}g_m(s_0-s,\lambda_0-\lambda)|^2\,|V_{g_n}f_2(t,\nu)|^2
\end{align*}
Defining
$$
w(s,\lambda) =\|f_2\|_2^{-2}\,\sum_{(t,\nu)
\in\Lambda_{(s,\lambda)}} |V_{g_n}f_2(t,\nu)|^2
\quad\text{and} \quad
\mu = \sum_{(s,\lambda)\in \Pi_1(\Lambda)}\, w(s,\lambda)
\,\delta_{(s,\lambda)}
$$
for $f_2\ne0$,
the above identity can be written as
$$
\int_{{\Bbb R}^m}\,|f_1(x)|^2\,dx= \sum_{(s,\lambda)\in\Pi_1(\Lambda)}\,
w(s,\lambda)\,
|V_{f_1}g_m(s_0-s,\lambda_0-\lambda)|^2 = \left(\mu\ast|V_{f_1}g_m|^2\right)(s_0,\lambda_0).
$$
On the other hand, letting $\check \chi_{[0,1)^{2m}}(s,\lambda)=\chi_{[0,1)^{2m}}(-s,-\lambda)$
and defining $C$ and $\Lambda(C)$ as above,
we have also
$$
\begin{aligned}
\left(\mu\ast \check\chi_{[0,1)^m}\right)(s_0,\lambda_0) &= \sum_{(s,\lambda)\in \Pi_1(\Lambda)}\,
w(s,\lambda)\,\chi_{[0,1)^{2m}}(s-s_0,\lambda-\lambda_0)\\
=&\sum_{(s,\lambda)\in \Pi_1(\Lambda)\cap C}\, w(s,\lambda)\\
&=\|f_2\|_2^{-2}\,\sum_{(t,\nu)\in\Lambda(C)}
 |V_{g_n}f_2(t,\nu)|^2 \leq 1,
\end{aligned}
$$
where the last inequality results from the orthogonality
of the system $(\chi_{[0,1]^n},\Lambda(C))$ proved earlier.
Since $(s_0,\lambda_0)$ is arbitrary in ${\Bbb R}^{2m}$
and
$$
\int_{{\Bbb R}^{2m}}|V_{f_1}g_m(s,\lambda)|^2\,ds\,d\lambda=\|f_1\|_2^{2},
$$
Theorem \ref{tiling Prop} can be used to deduce that $\mu\ast \check\chi_{[0,1)^m} =1$.
This shows that
$$
\sum_{(t,\nu)\in\Lambda(C)}
 |V_{g_n}f_2(t,\nu)|^2= \|f_2\|^2,
\quad f_2\in L^2({\Bbb R}^n),
$$
and thus that the  system $(\chi_{[0,1]^n},\Lambda(C))$ is complete, proving
our claim.
\end{pf}
\medskip

\noindent{\bf Proof of Theorem \ref{pseudo th}.}
Let ${\mathcal J} = \Pi_1(\Lambda)$ and, for any $(s,\lambda)\in{\mathcal J}$, define
$$
 \Lambda_{(s,\lambda)}=\{(t,\nu):(s,t,\lambda,\nu)\in \Lambda \}.
$$
If $(s_0,\lambda_0)\in{\mathcal J}$,
let $C = (s_0,\lambda_0)+[0,1)^{2m}$, and
$$
\Lambda(C): = \{(t,\nu)\in{\Bbb R}^{2n}: (s,t,\lambda,\nu)
\in \Lambda \ \mbox{and} \ (s,\lambda)\in C\}.
$$
Proposition  \ref{prop3.3} shows that the
 system $(\chi_{[0,1]^n},\Lambda(C))$
 forms a Gabor orthonormal basis. By assumption ${\mathcal J} +[0,1)^{2m}$
tiles ${\Bbb R}^{2m}$. Hence, $(s_0,\lambda_0)+[0,1)^{2m}$ contains exactly one point
in ${\mathcal J}$, i.e.~$(s_0,\lambda_0)$, and we have
$$
\Lambda(C)=
 \{(t,\nu):(s_0,t,\lambda_0,\nu)\in \Lambda \}
= \Lambda_{(s_0,\lambda_0)}.
$$
Therefore, we can write $\Lambda$ as
$$
\Lambda = \bigcup_{(s_0,\lambda_0)\in{\mathcal J}} \{(s_0,\lambda_0)\}\times \Lambda_{(s_0,\lambda_0)}.
$$
Our proof will be complete if we can show that ${\mathcal J}$ is a
Gabor orthonormal basis of $L^2({\Bbb R}^m)$.

\medskip

As ${\mathcal J}$ is a tiling set, by Proposition \ref{prop3.1}
it suffices to show that the inclusion
${\mathcal J}-{\mathcal J}\subset{\mathcal Z}
(V_{g_m}g_m)\cup\{0\}$ holds. Let $(s,\lambda)$ and $(s',\lambda')$
be distinct points in ${\mathcal J}$. As $\Lambda_{(s,\lambda)}+[0,1)^{2n}$
tiles ${\Bbb R}^{2n}$, so does $\Lambda_{(s,\lambda)}+[-1,0)^{2n}$,
and we can find
$(t,\nu)\in\Lambda_{(s,\lambda)}$ such that
$0\in (t,\nu)+[-1,0)^{2n}$, or, equivalently, with
$(t,\nu)\in[0,1)^{2n}$.
Similarly, we can find $(t',\nu') \in\Lambda_{(s',\lambda')}$
such that $(t',\nu')\in[0,1)^{2n}.$ Using the fact that
$(\chi_{[0,1]^d},\Lambda)$ is a Gabor
orthonormal basis of $L^2({\Bbb R}^{2d})$, we have
$$
(s,t,\lambda,\nu)-(s',t',\lambda',\nu')
\in{\mathcal Z}(V_{g_d}g_d).
$$
or, equivalently,
$$
V_{g_m}g_m(s-s',\lambda-\lambda')=0\quad
\text{or}\quad  V_{g_n}g_n(t-t',\nu-\nu')=0.
$$
Note that, since  $|t-t'|<1$ and $|\nu-\nu'|<1$,  $V_{g_n}g_n(t-t',
\nu-\nu')\neq0$. Hence $(s,\lambda)-(s',\lambda')
\in{\mathcal Z}(V_{g_m}g_m)$ as claimed.
\qquad$\Box$
\medskip

\section{Two-dimensional Gabor orthonormal bases}
In this section, our goal will be to classify all possible Gabor orthonormal basis
generated by the unit square on ${\Bbb R}^2$.

Given a fixed Gabor orthonormal
basis ${\mathcal G}(\chi_{[0,1]^2},\Lambda)$ and a set $A\subset {\Bbb R}^2$, we define
the sets
$$
\Gamma(A) = \{(\lambda_1,\lambda_2)\in{\Bbb R}^2 :
(t_1,t_2,\lambda_1,\lambda_2)\in\Lambda, \ (t_1,t_2)\in A\},
$$
 and, for any $(\lambda_1,\lambda_2)\in{\Bbb R}^2$ and any set $B\subset {\Bbb R}^2$, we let
 $$
T_{A}(\lambda_1,\lambda_2) = \{(t_1,t_2)\in{\Bbb R}^2:
(t_1,t_2,\lambda_1,\lambda_2)\in\Lambda,\,\, (t_1,t_2)\in A\}
$$
and
$$
T_{A}(B)=\{(t_1,t_2)\in{\Bbb R}^2:
(t_1,t_2,\lambda_1,\lambda_2)\in\Lambda,\,\, (t_1,t_2)\in A,\,\,(\lambda_1,\lambda_2)\in B\}.
$$

In particular, the set $T_{A}(\Gamma(A))$ collects all  the couples
$(t_1,t_2)\in A$ such that $(t_1,t_2,\lambda_1,\lambda_2)\in \Lambda$
for some $(\lambda_1,\lambda_2)\in{\Bbb R}^2$.

We say that a square is {\it half-open} if it is  a translate of one of the sets
$$
[0,1)^2,\quad (0,1]^2,\quad [0,1) \times(0,1]\quad\text{or}\quad(0,1]\times[0,1).
$$
Two measurable subsets of ${\Bbb R}^d$ will be called {\it essentially disjoint}
if their intersection has zero Lebesgue measure.
In the derivation below, we will make use of the  identity
$$
V_{g_2}g_2(t_1,t_2,\lambda_1,\lambda_2) = V_{g_1}g_1(t_1,\lambda_1) V_{g_1}g_1(t_2,\lambda_2),
\quad (t_1,t_2,\lambda_1,\lambda_2)\in {\Bbb R}^4,
$$
which implies, in particular, that
$$
V_{g_2}g_2(t_1,t_2,\lambda_1,\lambda_2)=0 \quad \iff
\quad V_{g_1}g_1(t_1,\lambda_1)=0\,\,\text{or}\,\,  V_{g_1}g_1(t_2,\lambda_2)=0.
$$
Moreover, using (\ref{ZVggd}), the zero set of $V_{g_2}g_2$ is given by
\begin{equation}\label{ZVgg2}
{\mathcal Z}(V_{g_2}g_2) = \{(t,\lambda): |t|_{\max}\geq 1\}\cup
\left(\bigcup_{i=1}^{2}\{(t,\nu): \lambda_i(1-|t_i|)\in{\Bbb Z}
\setminus\{0\})\}\right).
\end{equation}
This implies that if $|t|_{\max}<1$
and $(t,\lambda)\in{\mathcal Z}(V_{g_2}g_2)$,
then, there exists $i\in \{1,2\}$ and for some integer $m\ne 0$ such that
$$
|\lambda_i| = \frac{|m|}{1-|t_i|}\ge1.
$$
with a strict inequality if $t_i\ne 0$.
These properties will be used throughout this section.
\medskip

\begin{Lem}\label{lem4.1} Let  ${\mathcal G}(\chi_{[0,1]^2},\Lambda)$ be
a Gabor orthonormal basis for $L^2( {\Bbb R}^2)$
and let $C$ be a half-open square. Then,

\medskip

(i) $\Gamma(C)+[0,1]^2$ is a packing of ${\Bbb R}^2$.

\medskip

(ii)  If $(\lambda_1,\lambda_2)\in\Gamma(C)$, then
$T_{ C}(\lambda_1,\lambda_2) $ consists of one point.
\end{Lem}

\medskip

\begin{pf}
(i) Let $(\lambda_1,\lambda_2)$ and $(\lambda_1',\lambda_2')$ be distinct elements
of $\Gamma(C)$. By definition, we can find $(t_1,t_2)$ and $(t_1',t_2')$ in $C$ such
that  $(t_1,t_2,\lambda_1,\lambda_2), (t_1',t_2',\lambda_1',\lambda_2')\in\Lambda$.
We then have
$$
0 = V_{g_1}g_1(t_1-t_1',\lambda_1-\lambda_1')\,V_{g_1}g_1(t_2-t_2',\lambda_2-\lambda_2')
$$
If, without loss of generality, the first factor on the right-hand side of the previous equality
vanishes, the fact that $|t_1-t_1'|<1$ shows the existence of an integer $k> 0$ such that
$$
|\lambda_1-\lambda_1'| = k/(1-|t_1-t_1'|)\ge 1.
$$
Hence, the cubes
$(\lambda_1,\lambda_2)+[0,1]^2$ and $(\lambda_1',\lambda_2')+[0,1]^2$ are essentially disjoint.

\medskip

(ii) Suppose that $T_C(\lambda_1,\lambda_2)$ contains two distinct points $(t_1,t_2)$
and $(t_1',t_2')$. Then,
$$
0 = V_{g_1}g_1(t_1-t_1',0)\,V_{g_1}g_1(t_2-t_2',0).
$$
As $V_{g_1}g_1 (t,0)\neq 0$ for any $t$ with $|t|<1$, we must have $|t_1-t_1'|\geq 1$
or $|t_2-t_2'|\geq 1$, contradicting the fact that both $(t_1,t_2)$
and $(t_1',t_2')$ belong to $C$.
\end{pf}
In the following, we will denote by $\partial A$ the boundary of a set $A$.
The next result will be useful.
\medskip

\begin{Lem}\label{lem4.2-}
Under the hypotheses of the previous lemma, consider an element
$\lambda = (\lambda_1,\lambda_2)$ of ${\Gamma}(C)$ and  let $T_{C}(\lambda)
= \{(t_1,t_2)\}$. Then for any $x\in\partial(\lambda+[0,1]^2)$, we can find
$\lambda_x = (\lambda_{1,x},\lambda_{2,x})\in\Gamma(C)$ such that
$x\in\partial(\lambda_x+[0,1]^2)$.
Moreover, for any such $\lambda_x$,  letting $T_{C}(\lambda_x) = \{t_x\}$, where $t_x = (t_{1,x},t_{2,x})$,
we can find $i_0\in\{1,2\}$ such that $t_{i_0,x}=t_{i_0}$ and
$\lambda_{i_0,x} = \lambda_{i_0}+1$ or $\lambda_{i_0}-1$.
\end{Lem}

\begin{pf}
We can write $x = (\lambda_1+\epsilon_1,\lambda_2+\epsilon_2)$,
where $0\leq \epsilon_i\leq 1$, $i=1,2$ and $\epsilon_i\in\{0,1\}$ for at least
one index $i$. Let $a=(a_1,a_2)\in {\Bbb R}^2$ with $0< a_i<1$ for $i=1,2$ and
consider the point $(t_a,x) : = (t_1+a_1,t_2+a_2,\lambda_1+\epsilon_1,
\lambda_2+\epsilon_2)$ in ${\Bbb R}^{4}$. Since
$\Lambda+[0,1]^{4}$ is a tiling on ${\Bbb R}^{4}$ and  the point
$(t_a,x)$ is a point on the boundary of $(t,\lambda)+[0,1]^{4}$,
we can find some point $ (t_{x,a},\lambda_{x,a})\in \Lambda\setminus \{(t,\lambda)\}$
such that
$(t_a,x)\in(t_{x,a},\lambda_{x,a})+[0,1]^{4}$. Let $t_{x,a} =
(t_1',t_2')$ and $\lambda_{x,a} = (\lambda_1',\lambda_{2}')$. We have
\begin{equation}\label{eq4.0}
\left\{
  \begin{array}{ll}
    -a_i\leq t_i-t_i'\leq1-a_i , \\
-\epsilon_i\leq \lambda_i-\lambda_i'\leq 1-\epsilon_i,
  \end{array}
\right. \ i=1,2
\end{equation}
Using the orthogonality of the system ${\mathcal G}(\chi_{[0,1]^2},\Lambda)$,
 we can find $i_0\in\{1,2\}$ such
that  $V_{g_1}g_1(t_{i_0}-t_{i_0}',\lambda_{i_0}-\lambda_{i_0}')=0$.
Note that $t_{i_0}-t_{i_0}'\neq0$ would imply that $|\lambda_{i_0}-\lambda_{i_0}'|>1$
which is impossible from (\ref{eq4.0}). Hence, $t_{i_0}=t_{i_0}'$ and
$\lambda_{i_0}-\lambda_{i_0}'\ne 0$.

Moreover, as $V_{g_1}g_1(0,v)\neq 0$ if $|v|<1$, $V_{g_1}g_1(t_{i_0}-t_{i_0}',
\lambda_{i_0}-\lambda_{i_0}')=0$ can only occurs if $|\lambda_{i_0}-\lambda_{i_0}'|=1$.
This shows also that $\epsilon_{i_0}\in\{0,1\}$ in that case.
 This proves the last statement of our claim and the fact that
$x\in\partial(\lambda_{x,a}+[0,1]^2)$. The proof will be complete
if we can show that $\lambda_{x,a}\in\Gamma(C)$ for some choice of $a$.

\medskip

For simplicity, we consider the half-open square to be
${C} = [b_1,b_1+1)\times[b_2,b_2+1)$.  Our assertion will be true if the
point $t_{x,a} =
(t_1',t_2')$ constructed above satisfies
the inequalities $b_i\leq t_i'< b_{i}+1$ for $i=1,2$.
As $t_{i_0} = t_{i_0}'$, the inequalities clearly hold for $i=i_0$.
Suppose that the other index $j$ falls out of the range,
say $t_j'<b_j$ (The case $t_j'\geq b_{j}+1$ is similar).
We consider $(t_{a'},x)$ with $a_j' = t_j'+1-t_j+\delta$
for some small $\delta>0$. Note that, by (\ref{eq4.0}),
 we have $t_i+a_i-1\leq t_i'\leq t_i+a_i$ for $i=1,2$,
and, in particular,
$$
a_j'=t_j'+1-t_j+\delta\ge a_j+\delta>0.
$$
We have also $a_j' <1$. Indeed, the inequality
$t_j'-t_j+1+\delta \geq 1$ would imply that $t_j'+1 +\delta \geq 1+t_j.$
This is not possible, as $b_j \leq t_j<b_j+1$, so $1+t_j\geq b_j+1$. But
$t_j'<b_j$, so $t_j'+1<b_j+1$, so for $\delta$ small,
$$
t_j'+1+\delta<b_j+1 \leq 1+t_j
$$
which yields a contradiction.

Using the previous argument with $a'$ replacing $a$, we guarantee the
existence of $t_j''$ such that  $ t_j'+\delta = t_j+a_j'-1\leq t_j''
\leq t_j+a_j' =  t_j'+1+\delta$ and the associated point $(t_{a'},
\lambda_{x,a'}) =(t_1'',t_2'',\lambda_1'', \lambda_2'')$  in
$\Lambda$ with the property that $x\in\partial(\lambda_{x,a'}+[0,1]^2)$
for some index $i_0'$  such that
$|\lambda_{i_0'}-\lambda_{i_0'}''|=1$,
$t_{i_0'} =t_{i_0'}''$ and $\epsilon_{i_0'}\in\{0,1\}$.
We claim that $t_j'' = t_j'+1$.
Now, $(t_1',t_2',\lambda_1',\lambda_2')$ and
$(t_1'',t_2'',\lambda_1'', \lambda_2'')$ are in $\Lambda$.
The mutual orthogonality property implies that
$V_{g_1}g_1(t_i'-t_i'',\lambda_i'-\lambda_i'')=0$
for some $i=1,2$.

\medskip

Suppose that $x$ is not of the corner points of $\lambda+[0,1]^2$.
In that case, the index $i$ such that $\epsilon_i\in\{0,1\}$ is unique
and it follows that $i_0 =i_0'$. This implies in particular,
that $t_{i_0}' = t_{i_0}''$ ( as $t_{i_0}'=t_{i_0} = t_{i_0'} = t_{i_0'}'' = t_{i_0}''$).
Furthermore, the second set of inequalities in (\ref{eq4.0}) show that
$\lambda_{i_0}' =\lambda_{i_0}''=\lambda_{i_0}-1$ if $\epsilon_{i_0}=0$
and $\lambda_{i_0}' =\lambda_{i_0}''=\lambda_{i_0}+1$ if $\epsilon_{i_0}=1$.
We have thus $\lambda_{i_0}' = \lambda_{i_0}''$ in both cases.
We have thus
$$
V_{g_1}g_1(t_{i_0}'-t_{i_0}'',\lambda_{i_0}'-\lambda_{i_0}'')=
V_{g_1}g_1(0,0)=1.
$$
Therefore, the other index $j$ must satisfy
$V_{g_1}g_1(t_j'-t_j'',\lambda_j'-\lambda_j'')=0$.
The inequalities
$$-\epsilon_j\leq\lambda_j-\lambda_j'\leq 1-\epsilon_j
\quad\text{and}\quad-\epsilon_j\leq\lambda_j-\lambda_j''\leq 1-\epsilon_j
$$
yield $-1\leq \lambda_j'-\lambda_j''\leq 1$. However,
$\delta\leq t_j''-t_j'\leq 1+\delta$. The $V_{g_1}g_1$
would not be zero unless $t_j''\geq t_{j}'+1 (\geq b_j)$.
Hence, $t_j'+1\leq t_j''\leq t_j'+1+\delta$.
This forces that $t_j'' = t_j'+1$.  This completes
the proof for non-corner points. If $x$ is of the corner point,
as the square constructed for the non-corner will certainly
cover the corner point. Therefore, the proof is completed.
\end{pf}

\medskip

With the help of the previous two lemmas, the following  tiling result for
$\Gamma(C)$  follows immediately.

\begin{Cor}\label{Cor4.2}
Let $C$ be a half-open square. Then $\Gamma(C)+[0,1]^2$
is a tiling of ${\Bbb R}^2$.
\end{Cor}

\begin{pf}
It suffices to prove the following statement: suppose that ${\mathcal J}+[0,1]^2$ is non-empty
packing of  ${\Bbb R}^2$. If, for any $x\in\partial(t+[0,1]^2)$ where $t\in {\mathcal J}$, we can find
$t_x \in{\mathcal J}$ with $t_x\ne t$ such that
$x\in\partial(t_x+[0,1]^2)$, then  ${\mathcal J}+[0,1]^2$ is a tiling of  ${\Bbb R}^2$.
Indeed, by Lemma \ref{lem4.1}(i) and Lemma \ref{lem4.2-},
$\Gamma(C)+[0,1]^2$ is a packing of  ${\Bbb R}^2$ and satisfies the stated property. It is thus a tiling
of  ${\Bbb R}^2$.

\medskip

To prove the previous statement, we note that as ${\mathcal J}+[0,1]^2$ is packing, it is a closed set.
Suppose that ${\mathcal J}+[0,1]^2$ satisfies the property above and that
${\Bbb R}^d\setminus ({\mathcal J}+[0,1]^2)\neq \emptyset$.
Let $x\in \partial ({\mathcal J}+[0,1]^2)$ and assume that $x\in t+[0,1]^2$.
We can then find $t_x\in{\mathcal J}$ with $t_x\ne t$ such that
$x\in \partial (t_x+[0,1]^2)$.
Note that if $x$ were  not a corner point of either $t+[0,1]^2$ or $t_x+[0,1]^2$,
then $x$ would be  in the interior of ${\mathcal J}+[0,1]^2$. Hence,
$x$ must be a corner point of $t+[0,1]^2$ or $t_x+[0,1]^2$.
As the set of all the corner points of the squares in ${\mathcal J}+[0,1]^2$ is countable,
the Lebesgue measure of the open set ${\Bbb R}^d\setminus ({\mathcal J}+[0,1]^2)$ is zero
and ${\Bbb R}^d\setminus ({\mathcal J}+[0,1]^2)$ is thus empty, proving our claim.

\end{pf}


\medskip

\

\begin{Lem}\label{lem4.2} Let $C$ be a half-open square and suppose that
 $(\lambda_1,\lambda_2)\in{\Gamma}(C)$ with
$T_{C}(\lambda_1,\lambda_2) = \{(t_1,t_2)\}$.
Then all the sets  $T_{C}(\lambda_1',\lambda_2')$ with $(\lambda_1',\lambda_2')\in{\Gamma}(C)$
are either of the form
$\{(t_1,t_2+s)\}$ or $\{(t_1+s,t_2)\}$ for some  real $s$ with $|s|<1$
 depending on $(\lambda_1,\lambda_2)$.

\end{Lem}

\begin{pf} We first make the following remark. If
$(\alpha_1,\alpha_2),(\beta_1,\beta_2) \in{\Gamma}(C)$
are such that the two squares $(\alpha_1,\alpha_2)+[0,1]^2$
and  $(\beta_1,\beta_2)+[0,1]^2$ intersect each other
and also both intersect a third square
$(\gamma_1,\gamma_2)+[0,1]^2$ with $(\gamma_1,\gamma_2) \in{\Gamma}(C)$,
then, letting $T_{C}(\gamma_1,\gamma_2)=\{(r_1,r_2)$, we have
$$
T_{C}(\alpha_1,\alpha_2)=\{(r_1+a,r_2)\}\quad\text{and}\quad
T_{C}(\beta_1,\beta_2)=\{(r_1+b,r_2)\}
$$
or
$$
T_{C}(\alpha_1,\alpha_2)=\{(r_1,r_2+a)\quad\text{and}\quad T_{C}(\beta_1,\beta_2)=\{(r_1,r_2+b)\},
$$
for some real $a,b$.  Indeed, using Lemma \ref{lem4.2-}, we have
$T_{C}(\alpha_1,\alpha_2)=\{(r_1+a,r_2)\}$ or $\{(r_1,r_2+a)\}$
and  $T_{C}(\beta_1,\beta_2)=\{(r_1+b,r_2)$ or $\{(r_1,r_2+b)\}$.
Suppose, for example, that $T_{C}(\alpha_1,\alpha_2)=\{(r_1+a,r_2)\}$
and $T_{C}(\beta_1,\beta_2)=\{(r_1,r_2+b)$. Since the two squares intersect each other,
we must have $|\alpha_1-\beta_1|\leq 1$ and $|\alpha_2-\beta_2|\leq 1$.
The  orthogonality property also implies
 that either
$(a,\alpha_1-\beta_1)$ or $(-b,\alpha_2-\beta_2)$
is in the zero set of $V_{g_1}g_1$. But since we have $|a|, |b|<1$,
this would imply that $|\alpha_1-\beta_1|>1$ or $|\alpha_1-\beta_2|>1$,
which cannot happen.
As $\Gamma(C)+[0,1]^2$
is a tiling of ${\Bbb R}^2$, for any square
$(\sigma_1,\sigma_2)+[0,1]^2$ intersecting
the square  $(\lambda_1,\lambda_2)+[0,1]^2$ and with  $(\sigma_1,\sigma_2)\in\Gamma(C)$,
we can find another square  $(\delta_1,\delta_2)+[0,1]^2$, with
$(\delta_1,\delta_2)\in\Gamma(C)$  and with $(\delta_1,\delta_2)+[0,1]^2$
intersecting both squares $(\sigma_1,\sigma_2)+[0,1]^2$ and $(\lambda_1,\lambda_2)+[0,1]^2$.
By the previous remark, the conclusion of the lemma holds for all the squares
that neighbour  the square $(\lambda_1,\lambda_2)+[0,1]^2$. Replacing this original
square by one of the neighbouring squares and continuing this process,
we obtain the conclusion of the lemma
for all the squares in the tiling  $\Gamma(C)+[0,1]^2$ by an induction argument.
This proves our claim.
\end{pf}

\medskip


Suppose that the system ${\mathcal G}(\chi_{[0,1]^2},\Lambda)$ gives rise to a
non-standard Gabor orthonormal basis of $L^2({\Bbb R}^2)$.
Then, some of the squares will have overlaps and, without loss of
generality, we can assume that
$$
|[0,1]^2\cap [0,1]^2+(t_1,t_2)|>0
$$
for some $(t_1,t_2)$ in the translation  component of $\Lambda$.
\medskip

\begin{Lem}\label{Lem4.3}
 If $(0,0,0,0)\in\Lambda$, then the sets
$T_{[0,1)^2}(\lambda_1,\lambda_2)$ where
$(\lambda_1,\lambda_2)\in\Gamma([0,1)^2)$ are either all of the form
$ \{(t,0)\}$ or all w of the form $\{(0,t)\}$ with some $t$ (depending on
$(\lambda_1,\lambda_2)$) with $|t|<1$.
In the first case, if there exists some $(\lambda_1,\lambda_2)\in\Gamma([0,1)^2)$
with $T_{[0,1)^2}(\lambda_1,\lambda_2)=(t,0)$ and
$t\neq0$, then
\begin{equation}\label{eqGamma}
\Gamma([0,1)^2) = \bigcup_{k\in{\Bbb Z}}({\Bbb Z}+\mu_{k,0})\times\{k\}
\end{equation}
for some $0\leq\mu_{k,0}<1$. Moreover, we can find $0\leq t_k<1$ such that
\begin{equation}\label{eqT[0,1]}
T_{[0,1)^2}(({\Bbb Z}+\mu_{k,0})\times\{k\}) = \{(t_k,0)\},\quad k\in {\Bbb Z},
\end{equation}
and
\begin{equation}\label{eqLambda01}
\Lambda\cap ([0,1)^2\times {\Bbb R}^2) = \{(t_k,0, j+\mu_{k,0}, k): j,k\in{\Bbb Z}\}.
\end{equation}
(In the second case,
$\Gamma([0,1)^2) = \bigcup_{k\in{\Bbb Z}}\{k\}\times({\Bbb Z}+\mu_{k,0})$
and $T_{[0,1)^2}(\{k\}\times({\Bbb Z}+\mu_{k,0})) =
\{(0,t_k)\}$, $\Lambda\cap ([0,1)^2\times {\Bbb R}^2) = \{(0,t_k,k,j+\mu_{k,0}):j,k\in{\Bbb Z}\})$.
\end{Lem}

\begin{pf}
If $\lambda = (0,0)$, we have $T_{[0,1)^2}(\lambda) = \{(0,0)\}$
as  $(0,0,0,0)\in \Lambda$.
By Lemma \ref{lem4.2}, any $(\lambda_1,\lambda_2)\in\Gamma([0,1)^2)$
with the square $(\lambda_1,\lambda_2)+[0,1]^2$ intersecting $[0,1]^2$ on
the $\lambda_1,\lambda_2$-plane satisfies  $T_{[0,1)^2}(\lambda_1,\lambda_2)=\{(t,0)\}$
 or $T_{[0,1)^2}(\lambda_1,\lambda_2)\{(0,t)\}$ with $|t|<1$.
 Without loss of generality,
we assume that the first case holds. As $\Gamma([0,1)^2)+[0,1]^2$
is a tiling of ${\Bbb R}^2$, for any square
$C = (\lambda_1,\lambda_2)+[0,1]^2$, with  $(\lambda_1,\lambda_2)\in\Gamma([0,1)^2)$,
 we can find squares $C_{i} = (\lambda_{1,i},\lambda_{2,i})+[0,1]^2$
for $i=0,\dots,k$ with $(\lambda_{1,i},\lambda_{2,i})\in\Gamma([0,1)^2)$
and such that
$C_0=[0,1]^2$, $C_k = C$,
and with $C_i$ and $C_{i+1}$ touching each other for all $i=0,\dots, k-1$.

We have $T_{[0,1)^2}(\lambda_{1,1},\lambda_{2,1}) = \{(t_1,0)\}$
for some number $t_1$ with $|t_1|<1$.
Since $C_2$ and $C_0$ both intersect $C_1$,
$T_{[0,1)^2}(\lambda_{1,2},\lambda_{2,2}) = \{(t_2,0)\}$
by Lemma \ref{lem4.2} again.  Inductively, we have
$T_{[0,1)^2}(\lambda_{1,i},\lambda_{2,i}) = \{(t_i,0)\}$, $i=1,\dots,k$,
which proves the first part.

\medskip

Consider the case where, for any $(\lambda_1,\lambda_2)\in\Gamma([0,1)^2)$,
there exists a number $t=t(\lambda_1,\lambda_2)$ such that
$T_{[0,1)^2}(\lambda_{1},\lambda_{2})
= \{(t,0)\}$ and assume that  $t(\lambda_1,\lambda_2)\ne 0$ for at least
one couple  $(\lambda_1,\lambda_2)\in\Gamma([0,1)^2)$.
Suppose that $\Gamma([0,1)^2)$ is not of the form in (\ref{eqGamma}).
By  Corollary \ref{Cor4.2} and Proposition \ref{prop3.2}, we must have
$\Gamma([0,1)^2) = \bigcup_{k\in{\Bbb Z}}\{k\}\times({\Bbb Z}+a_k)$
with $0\le a_k<1$ and at least one  $a_k\neq0$.
Consider  the distinct  points
$$
(t,0,k,a_k+j)\quad \text{and}\quad (t',0,k-1,a_{k-1}+j),\ \ \text{both in}\,\,\Lambda.
$$
We must have that either $(t-t',1)\in
{\mathcal Z}(V_{g_1}g_1)$ or $(0,a_k-a_{k-1})\in{\mathcal Z}(V_{g_1}g_1)$.
However, since $|a_k-a_{k-1}|<1$, the second case is impossible.
This means that $(t-t',1)\in {\mathcal Z}(V_{g_1}g_1)$ which is
possible only if $t=t'$. Therefore the fact that $(t,0,k,a_k+j)\in \Lambda$ implies that $t=t_j$
for some real $t_j$. We know prove by induction on $|j|$ that $t_j=0$ for all $j\in \mathbb{Z}$.
The case $j=0$ is clear as $(0,0,0,0)\in \Lambda$ by assumption.
If our claim is true for all $|j|\le J$ where $J\ge 0$, chose $k\in \mathbb{Z}$ such that
$a_{k+1}\ne 0$ and $a_k=0$ if such $k$ exists.
Suppose first that $j>0$. There exist thus $t\in [0,1)$ such that
$$
(t_{j+1},0,k,j+1)\quad \text{and}\quad (0,0,k+1,a_{k+1}+j)\ \ \text{both belong  to}\,\,\Lambda.
$$
This implies that either $(t,-1)\in {\mathcal Z}(V_{g_1}g_1)$
or $(0,a_{k+1}-1)\in {\mathcal Z}(V_{g_1}g_1)$. This last case is impossible and the first one
is only possible if $t=0$, showing that $t_{j+1}=0$.
Similarly by considering the points
$$
(t_{j-1},0,k+1,a_{k+1}+j-1)\quad \text{and}\quad (0,0,k,j)\ \ \text{which both belong  to}\,\,\Lambda.
$$
we can conclude that  $t_{j-1}=0$ for $j<0$. If $k$ as above does not exist, there exists
chose $k'\in {\Bbb Z}$ such that
$a_{k'-1}\ne 0$ and $a_{k'}=0$. By considering the points
$$
(t_{j+1},0,k',j+1)\quad \text{and}\quad (0,0,k'-1, a_{k'-1}+j)\quad \text{if}\,\,j>0
$$
and the points
$$
(t_{j-1},0,k'-1,a_{k'-1}+j-1)\quad \text{and}\quad (0,0,k', j)\quad \text{if}\,\,j<0
$$
which all belong to $\Lambda$, we conclude that $t_j=0$ if $|j|=J+1$.
This proves  (\ref{eqGamma}).

\medskip
If we are in the first case, i.e.
$$
\Gamma([0,1)^2) = \bigcup_{k\in{\Bbb Z}}({\Bbb Z}+\mu_{k,0})\times\{k\},
$$
let $m,m'$ be distinct integers. We have then
$$
T_{[0,1)^2}(m+\mu_{n,0},n) = \{(t_{m},0)\}\quad \text{and}
\quad
T_{[0,1)^2}(m'+\mu_{n,0},n) = \{(t_{m'},0)\}
$$
which implies that $V_{g_1}g_1(t_m-t_{m'},m-m') =0$ or $V_{g_1}g_1(0,0)=0$.
The second case is clearly impossible while  the first one is possible only
when $t_m=t_m'$. This shows (\ref{eqT[0,1]}) and (\ref{eqLambda01})
follows immediately from (\ref{eqGamma}) and (\ref{eqT[0,1]}).
\end{pf}

 Note that Lemma \ref{Lem4.3} implies that $\Gamma ([0,1)^2) =
\Gamma(\{(x,0):0\leq x<1\})$ and $\Gamma ((0,1)^2)=\emptyset$ if $(0,0,0,0)\in \Lambda$.

\medskip

\begin{Lem}\label{Lem4.4}
Under the assumptions of Lemma \ref{Lem4.3},
suppose that there exists $(\lambda_1,\lambda_2)\in\Gamma([0,1)^2)$
with $T_{[0,1)^2}(\lambda_1,\lambda_2)=(t,0)$ and
$t\neq0$.
Then we can find numbers $t_k$ with $0\le t_k<1$ and $\mu_{k,m},\,\, k,m\in {\Bbb Z}$,
with $0\leq \mu_{k,m}<1$,  such that
$$
\Lambda\cap ({\Bbb R}\times[0,1)\times {\Bbb R}^2) = \{(m+t_k,0, j+\mu_{k,m}, k): j,k,m \in{\Bbb Z}\}
$$
\end{Lem}

\medskip

\begin{pf}
By the result of  Lemma \ref{Lem4.3}, we have the identities
(\ref{eqT[0,1]}) and (\ref{eqLambda01}).
Let $T=\{t_k,\,\,k\in {\Bbb Z}\}\subset [0,1)$ where  $t_k$, $k\in{\Bbb Z}$,
are the numbers appearing in (\ref{eqT[0,1]}).
Let $s_1, s_2\in T$ with $s_1< s_2$. Consider the
half-open squares $C= (s_1,0)+[0,1)^2$ and $C' = (s_1,0)+\left((0,1]\times[0,1)\right)$.
Then we know that $\Gamma(C)+[0,1]^2$ and $\Gamma(C')+[0,1]^2$ both tile ${\Bbb R}^2$.
Let $P_0 = \{(s_1,y): 0\leq y<1\}$ and
$P_1 =\{(s_1+1,y):0\leq y< 1\} $.
Note that
$\Gamma(P_0) = \Gamma(\{(s_1,0)\})$.
Moreover,
$$
\Gamma(C) = \Gamma(P_0)\cup \Gamma(C\setminus P_0), \ \Gamma(C')
= \Gamma(C'\setminus P_1)\cup  \Gamma(P_1)
$$
and since $C\setminus P_0 = C'\setminus P_1$, $\Gamma(P_0)=\Gamma(P_1)$.
We have
$$
T_{C'}(\Gamma(P_1))\subset \{(s_1+1,y),\,\,0\le y<1\}
$$
but since $(s_2,0)\in C'$,
we must have $T_{C'}(\Gamma(P_1))=(s_1+1,0)$ by Lemma \ref{lem4.2}.
Since
$$
\Gamma(P_0) = \{(j+\mu_{k,0},k):\,\,j,k\in {\Bbb Z},\,\, t_k=s_1\}
$$
and
$\pi_2(\Gamma(P_0)) = \pi_2(\Gamma(P_1))$,
where $\pi_2$ is the projection to the second coordinate,
we have
$$
\Gamma(\{(1+s_1,0)\}) = \Gamma(P_1) = \{(j+\mu_{k,1},k):\,\,j,k\in {\Bbb Z},\,\,  t_k=s_1\}.
$$
for some constants $ \mu_{k,1}$ with $0\leq \mu_{k,1}<1$
using Proposition \ref{prop3.2}. Applying this argument to  $s_1=0$ and $s_2=t$,
we obtain that
$$
\Lambda\cap\left(\{1\}\times [0,1)\times {\Bbb R}^2\right)
=\{(j+\mu_{k,1},k):\,\,j,k\in {\Bbb Z},\,\,  t_k=0\}.
$$
Similar arguments applied to $s_1=s$ and $s_2=1$ show that, for any $s\in T$, we have
$$
\Lambda\cap\left(\{s+1\}\times [0,1)\times {\Bbb R}^2\right)
=\{(j+\mu_{k,1},k):\,\,j,k\in {\Bbb Z},\,\,  t_k=s\}.
$$
and that $\Lambda\cap\left(\{s+1\}\times [0,1)\times {\Bbb R}^2\right)$ is empty
if $s\in [0,1)\setminus T$. The same idea can also be used to show the existence
of constants $ \mu_{k,-1}$ with $0\leq \mu_{k,1}<1$ such that
$$
\Lambda\cap\left(\{s-1\}\times [0,1)\times {\Bbb R}^2\right)
 =\begin{cases}
\{(j+\mu_{k,-1},k):\,\,j,k\in {\Bbb Z},\,\,  t_k=s\},&  s\in T,\\
\qquad \emptyset,& s\in [0,1)\setminus T .
\end{cases}
$$
and, more generally using induction, that, for any $m\in {\Bbb Z}$,
we can find constants $\mu_{k,m}$ with $0\leq \mu_{k,m}<1$ such that
$$
\Lambda\cap\left(\{s+m\}\times [0,1)\times {\Bbb R}^2\right)
 =\begin{cases}
\{(j+\mu_{k,m},k):\,\,j,k\in {\Bbb Z},\,\,  t_k=s\},&  s\in T,\\
\qquad \emptyset,& s\in [0,1)\setminus T .
\end{cases}
$$
This proves our claim.

\end{pf}

%

\medskip

We can now complete the proof of the main result of
this section which gives a characterization for the subsets $\Lambda$ of ${\Bbb R}^4$
with the property that the associated set
 of time-frequency shifts applied to the window $\chi_{[0,1]^2}$
yields an orthonormal basis for $L^2({\Bbb R}^2)$.

\noindent{\bf Proof of Theorem \ref{th0.3}.} It follows from
Lemma \ref{lem4.2} that either all
$T_{[0,1)^2}(\lambda_1,\lambda_2)$,
$(\lambda_1,\lambda_2)\in\Gamma([0,1)^2)$ are either of the
form $ \{(t,0)\}$ or all are of the form $\{(0,t)\}$ with
some $t\neq 0$. In the first case, we deduce from Lemma \ref{Lem4.4}
that
$$
\Lambda\cap ({\Bbb R}\times[0,1)\times {\Bbb R}^2) = \{(m+t_k,0, j+\mu_{k,m}, k): j,k,m \in{\Bbb Z}\}
$$
for certain numbers $t_k$ and $\mu_{k,m}$ in the interval $[0,1)$.
We now show that $\Lambda$ will be of the first   of the two possible
forms given in the theorem.
(Similarly, the second form follows from the second case of Lemma \ref{Lem4.4}).

\medskip

Letting $C = [0,1)^2$ and $C' =[0,1)\times(0,1]$,
we note that both $\Gamma(C)+[0,1]^2$ and $\Gamma(C')+[0,1]^2$ tile
${\Bbb R}^2$ but $\Gamma ((0,1)^2)$ is empty. Hence,
$\Gamma(C') = \Gamma(\{(x,1): 0\leq x<1\})$.
It means that any set $T_{C'}(\lambda_1,\lambda_2)$
with $(\lambda_1,\lambda_2)\in\Gamma(C')$ is of the
form $\{(t,1)\}$ for some $t=t(\lambda_1,\lambda_2)$ with  $0\leq t<1$.
We now have two
possible cases: either the cardinality of $T_{C'}(\Gamma(C')$ is larger than one
or equal to one.
In the first case, we can find two distinct elements of
$T_{C'}(\Gamma(C'))$ and we can then replicate the
proof of Lemma \ref{Lem4.4} to obtain that
$$
\Lambda\cap ({\Bbb R}\times[1,2)\times {\Bbb R}^2)  = \{(m+t_k,1, j+\mu_{k,m,1}, k): j,k\in{\Bbb Z}\}.
$$
In the other case,  $T_{C'}(\Gamma(C'))=\{(t_1,1)\}$
for some $t_1$ with $0\le t_1<1$. If we translate $C'$ horizontally and use the same argument
as in the proof of Lemma \ref{Lem4.4}, we see that
$$
\Lambda\cap ({\Bbb R}\times[1,2)\times {\Bbb R}^2)  = \{(m+t_1,1)\}\times\Lambda_{m,1},
$$
where $\Lambda_{m,1}$ is a spectrum for the unit square $[0,1]^2$.
This last property is equivalent to $\Lambda_{m,1}+[0,1]^2$ being a tiling of ${\Bbb R}^2$
by the result in \cite{[LRW]}.

We can them prove  the theorem inductively by translating the square $C'$ in the vertical direction
using integer steps.
\qquad$\Box$
\medskip

\medskip

\end{document}